\documentclass{amsart}
\pdfoutput=1
\usepackage{graphicx}
\usepackage{subfigure}
\usepackage{url}
\usepackage{amsmath}
\usepackage{amssymb}
\usepackage{amsopn}
\usepackage{amssymb}
\usepackage{amstext}
\usepackage{amsthm}
\usepackage{float}
\usepackage{setspace}
\usepackage[T1]{fontenc}

\def \fig_size{ 4.5 }

\newtheorem{theorem}{Theorem}

\newcommand{\rank}{\operatorname{Rank}}

\newcommand{\sinc}{\operatorname{sinc}}

\let\oldsum\sum
\renewcommand{\sum}{\displaystyle\oldsum}
\let\oldmin\min
\renewcommand{\min}{\displaystyle\oldmin}
\let\oldmax\max
\renewcommand{\max}{\displaystyle\oldmax}
\let\oldinf\inf
\renewcommand{\inf}{\displaystyle\oldinf}
\let\oldsup\sup
\renewcommand{\sup}{\displaystyle\oldsup}

\graphicspath{ {./figures/} }

\begin{document}

\title{Robust Anomaly Detection Using Semidefinite Programming}

\author{J. A. Lopez}
\thanks{This work was supported in part by NSF grants IIS--1318145 and ECCS--1404163; AFOSR grant FA9550-12-1-0271the Alert DHS Center of Excellence under Award Number 2008-ST-061-ED0001.}
\author{O. Camps}
\author{M. Sznaier}
\thanks{The authors are with the Department of Electrical \& Computer Engineering, Northeastern University, MA 02115, USA. E-mails: {\tt lopez.jo@husky.neu.edu},{\tt \{camps,msznaier\}@coe.neu.edu }}


\begin{abstract} This paper presents a new approach, based on polynomial optimization and the method of moments, to the problem of anomaly detection.  The proposed technique only requires information about the statistical moments of the normal-state distribution of the features of interest and compares favorably with existing approaches (such as Parzen windows and 1-class SVM). In addition, it provides a succinct description of the normal state. Thus, it leads to a substantial  simplification of the the anomaly detection problem when working with higher dimensional datasets.   
 
\end{abstract}
\date{\today}

\maketitle

\section{Introduction}
Detecting anomalies in data is a task that engineers and scientists of every discipline encounter on a fairly regular basis. Common approaches to identifying anomalous observations include: estimating the probability density function of the ``normal'' state using density-estimation methods, computing the Mahalanobis distance between a point and a sample distribution, using machine-learning methods to learn the normal state and perform 1-class classification, using ``whiteness'' tests to detect when the differences between new samples and predictions (of linear models) become statistically colored, etc.  However,  a disadvantage encountered while estimating the distribution of the normal state  is that one never knows if the available samples are representative of \emph{all} the normal behaviors. This is especially true when there is a dearth of data with which to construct a reliable probability density estimate. In this paper, we address this challenge through the use of recent results from the method of moments and polynomial optimization.
The main idea is to use information about the statistical moments of the normal-state distribution of a given set of features  to compute an upper-bound 
 on the probability of a given observation of a
  random variable with this distribution.  The observation is then deemed to be anomalous when this bound falls below a given threshold.  While in principle computing this bound is a challenging infinite-dimensional problem, as shown in the sequel, the use of concepts from the theory of moments allows for recasting it into a computationally tractable finite dimensional optimization.

The generalized moment problem (GMP), defined below, has tremendous modeling power and has found use in a wide range of applications in areas  such as probability, financial economics, engineering, and operations research \cite{lasserre_text}. However, its wider application has been held back because, in its fullest generality, the GMP is intractable \cite{lasserre_sdp_approach_to_gpm}. On the other hand, recent advances in algebraic geometry have rendered the GMP tractable when the problem instance is described using polynomials \cite{lasserre_2001,lasserre_text}. 

Formally, the GMP is given by:
\begin{align}
\label{eq:gmp}
\begin{array}{cc}
\text{GMP:} & \rho_{\text{mom}} = \sup_{\mu\in\mathcal{M}(\mathbb{K})_+}\int_\mathbb{K}f_0 d\mu\\
\text{subject to } & \int_\mathbb{K} f_jd\mu \leqq \gamma_j, \; j=1,\dots,m
\end{array}
\end{align}
where $\mathbb{K}$ is a Borel subset of $\mathbb{R}^n$, $\mathcal{M}(\mathbb{K})$ is the space of finite signed Borel measures on $\mathbb{K}$, $\mathcal{M}(\mathbb{K})_+$ is the positive cone containing finite Borel measures $\mu$ on $\mathbb{K}$, $\lbrace\gamma_j\,|\,j=1,\dots,m\rbrace$ is a set of real numbers, and $f_j:\mathbb{K}\to\mathbb{R}$ are integrable with respect to all $\mu\in\mathcal{M}(\mathbb{K})_+$ for every $j=0,\dots,m$. The $\leqq$ stands for either an inequality or equality constraint.

Lasserre \cite{lasserre_2001,lasserre_text,lasserre_sdp_approach_to_gpm,lasserre_moms_and_sos_for_poly_optzn_and_related} showed that semidefinite programming (SDP) can be used to obtain an (often) finite sequence of relaxations which can be solved efficiently and which converge from above to the optimal solution when the problem has polynomial data; thereby giving us a tool with which to attempt to solve very difficult (i.e., non-convex) medium-sized problems, given the present state of SDP solvers. The monotonic convergence means that the $i$-th relaxation is useful, even when optimality has not been achieved, by providing the user with an upper (i.e. optimistic) bound $\rho_i$ to $\rho_{\text{mom}}$ in (\ref{eq:gmp}).

In the sequel, we will go over the preliminaries of the moments approach, present the aforementioned Lasserre relaxations, and demonstrate how to use them to detect anomalous data samples. 

\section{Preliminaries}
For a real symmetric matrix $A$, $A\ge 0$ means $A$ is positive semidefinite. Let $d\in\mathbb{N}\triangleq\{0,1,2,\dots\}$, $s(d)=\binom{n+d}{n}$, and let $\nu\in\mathbb{R}^{s(d)}$ be the column vector containing the canonical polynomial basis monomials
\begin{align}
\label{eq:polynomial_basis}
&\nu(x) = \left( 1,x_1,x_2,\dots, x_n, x_1^2, x_1 x_2,\dots,x_1 x_n,x_2 x_3,\dots,x_n^2,\dots,x_1^d,\dots,,x_n^d \right)^T
\end{align} 
for polynomials of the form $p(x):\mathbb{R}^n\to\mathbb{R}$ with dimension $s(d)$ and degree at most $d$. $\mathbb{R}[x]$ denotes the set of polynomials in $x \in \mathbb{R}^n = (x_1,x_2, \dots, x_n)$ with real coefficients. 

With $\alpha\in\mathbb{N}^n = (\alpha_1,\dots,\alpha_n)$, let  $x^\alpha = (x_1^{\alpha_1},\dots, x_n^{\alpha_n})$ and $y=\lbrace y_\alpha \rbrace$ denote a finite sequence of real variables. For a polynomial $p=\sum_\alpha p_\alpha x^\alpha$, with coefficients $p_\alpha$, let $L_y(p)$ denote the linear map
\begin{equation}
\label{eq:ly}
L_y(p) = \sum_\alpha p_\alpha y_\alpha
\end{equation}
which associates $p$ with the moment variables, $y_\alpha$.
\subsection{Moment and Localizing Matrices}
The $d$-th order moment matrix $M_d(y)$ is indexed by up-to order-$d$ monomial coefficients $\alpha,\beta\in\mathbb{N}^n$. That is, each entry of the $s(d)\times s(d)$ matrix is given by  
\begin{align}
M_d(y)(\alpha,\beta) = y_{\alpha+\beta},\; \alpha,\beta\in\mathbb{N}^n,\; |\alpha|,|\beta|\le d
\end{align}
where $|\alpha|=\sum_{i=1}^n \alpha_i$. 
For example, in the case $d=2$ and $n=2$, if $\alpha=(0,1)$ and $\beta=(2,0)$, then $M_2(y)(\alpha,\beta) = y_{21}$. Usually the $(\alpha,\beta)$ is suppressed (since it is understood that the first row and column of $M_d(y)$ contains the elements of $L_y(v)$). The entire $M_2(y)$ matrix is given below for illustration.
\begin{align}
M_2(y) &= \left[ 
\begin{array}{cccccc}
y_{00} & y_{10} & y_{01} & y_{20} & y_{11} & y_{02}\\
y_{10} & y_{20} & y_{11} & y_{30} & y_{21} & y_{12}\\
y_{01} & y_{11} & y_{02} & y_{21} & y_{12} & y_{03}\\
y_{20} & y_{30} & y_{21} & y_{40} & y_{31} & y_{22}\\
y_{11} & y_{21} & y_{12} & y_{31} & y_{22} & y_{13}\\
y_{02} & y_{12} & y_{03} & y_{22} & y_{13} & y_{04}
\end{array}
\right]
\end{align}
Thus, it is clear that $M_d(y)=L_y\left(\nu(x)\nu(x)^T\right)$. 

The classical problem of moments is concerned with determining if, for a given moment sequence $\lbrace y_\alpha\rbrace$, there is a measure $\mu$ so that $y_\alpha=\int x^\alpha d\mu$ for each $\alpha$. If so, we say $y$ has a \emph{representing} measure which is called \emph{determinate} if it is unique. Loosely speaking,  checking existence of  such a measure involves testing if $M_d(y)$ is positive semidefinite. Readers interested in the technical details are referred to \cite{lasserre_text} and references therein.

For a given polynomial $g\in\mathbb{R}[x]$, the \emph{localizing} matrix $M_d(g y)$ is given by
\begin{equation}
M_d(g y) = L_y(g(x)\nu(x)\nu(x)^T) 
\end{equation}
This is equivalent to shifting $M_d(y)$ by the monomials of $g$. For example, with $g(x)=a-x_1^2-x_2^2$, we have
\begin{align}
M_2(gy) &= a\left[ 
\begin{array}{cccccc}
y_{00} & y_{10} & y_{01} & y_{20} & y_{11} & y_{02}\\
y_{10} & y_{20} & y_{11} & y_{30} & y_{21} & y_{12}\\
y_{01} & y_{11} & y_{02} & y_{21} & y_{12} & y_{03}\\
y_{20} & y_{30} & y_{21} & y_{40} & y_{31} & y_{22}\\
y_{11} & y_{21} & y_{12} & y_{31} & y_{22} & y_{13}\\
y_{02} & y_{12} & y_{03} & y_{22} & y_{13} & y_{04}
\end{array}
\right]-
\left[
\begin{array}{cccccc}
y_{20} & y_{30} & y_{21} & y_{40} & y_{31} & y_{22}\\
y_{30} & y_{40} & y_{31} & y_{50} & y_{41} & y_{32}\\
y_{21} & y_{31} & y_{22} & y_{41} & y_{32} & y_{23}\\
y_{40} & y_{50} & y_{41} & y_{60} & y_{51} & y_{42}\\
y_{31} & y_{41} & y_{32} & y_{51} & y_{42} & y_{33}\\
y_{22} & y_{32} & y_{23} & y_{42} & y_{33} & y_{24}
\end{array}
\right]-\notag\\
&  \left[ 
\begin{array}{cccccc}
y_{02} & y_{12} & y_{03} & y_{22} & y_{13} & y_{04}\\
y_{12} & y_{22} & y_{13} & y_{32} & y_{23} & y_{14}\\
y_{03} & y_{13} & y_{04} & y_{23} & y_{14} & y_{05}\\
y_{22} & y_{32} & y_{23} & y_{42} & y_{33} & y_{24}\\
y_{13} & y_{23} & y_{14} & y_{33} & y_{24} & y_{15}\\
y_{04} & y_{14} & y_{05} & y_{24} & y_{15} & y_{06}
\end{array}
\right]
\end{align}
Localizing matrices are used to specify support constraints on $\mu$, as described in the next section.
\subsection{Upper-bounds Over Semi-algebraic Sets} The material in this subsection comes from \cite{lasserre_text} and \cite{lasserre_bounds_on_measures_with_moms} which discuss the problem of finding bounds on the probability that some random variable $x$ belongs to a set $S\subseteq\mathbb{R}^n$, given some of its moments $\gamma_\alpha$. 

Suppose $S$ is of the form
\begin{equation}
\label{eq:semialgebraic_set}
S = \left\lbrace x\in\mathbb{R}^n | f_j(x)\ge 0,\; j=1,\dots,m\right\rbrace
\end{equation}
where $f_j$ are given polynomial constraints. In terms of the GMP, this problem can be expressed as

\begin{align}
\label{eq:gmp_upper_bound}
\rho_{\text{mom}} &= \sup_{\mu\in\mathcal{M}(\mathbb{K})_+}\int_{\mathbb{R}^n} 1_S d\mu\\
\text{s.t. } & \int_{\mathbb{R}^n} x^\alpha d\mu = \gamma_j, \; j=1,\dots,m\notag
\end{align}
where $1_S$ is the indicator function of the set $S$. To solve this problem using SDP methods, Lasserre provided the following relaxations

\begin{align}
\mathbb{P}_S\mapsto
\begin{cases}
&\rho_i = \sup_{y,z} \; y_0\\
\text{s.t. } & y_\alpha + z_\alpha = \gamma_\alpha,\;\alpha\in\Gamma\\
& M_i(y) \ge 0\\
& M_i(z) \ge 0\\
& M_{i-v_j}(f_jy) \ge 0,\; j=1,\dots,m
\end{cases}
\label{eq:sdp_relaxations}
\end{align}
where $v_j = \lceil \deg(f_j)/2 \rceil$ and $\Gamma$ is the set of indices that correspond to the known moment information.

The following result assures us that $\rho_i$ approaches the optimal value $\rho_{\text{mom}}$ from above and tells us when optimality has been reached. This result is included for completeness as we will generally not require an optimal solution. This is an advantage because we do not have to worry about extracting optimizers\footnote{There is (currently) no reliable method of extracting optimizers from a moment matrix.} from a moment matrix that has satisfied the rank condition in Theorem \ref{theorem:main_theorem}.

\begin{theorem}[ ] \cite{lasserre_text} Let $\rho_i$ be the optimal value of the semidefinite program $\mathbb{P}_S$. Then:
\begin{description}
\item[(a)] For every $i>v$, $\rho_i\ge \rho_{\text{mom}}$ and moreover, $\rho_i\downarrow\rho_{\text{mom}}$ as $i\to\infty$.
\item[(b)] If $\rho_i$ is attained at an optimal solution $(y,z)$ which satisfies
\begin{center}
\begin{equation}
\left\lbrace 
\begin{array}{c}
\rank M_i(y)=\rank M_{i-v}(y)\\
\rank M_i(z)=\rank M_{i-1}(z)\\
\end{array}
\right.
\end{equation}
\end{center}
then $\rho_i=\rho_{\text{mom}}$.
\end{description}
\label{theorem:main_theorem}
\end{theorem}

\subsection{Lower Bounds} To compute lower bounds, one can use the complement of $S$, instead of $S$, in $\mathbb{P}_S$ and subtract $\rho_i$ from 1 \cite{lasserre_text}. However, in the approach presented here, the lower bounds do not provide much information because the sets $S$ that we will work with are small compared to the range of the data, which means the lower bound will usually be zero.

\section{Robust Anomaly Detection}
One approach to anomaly detection is to estimate the probability density function (PDF) and select a threshold below which incoming samples are classified as anomalous. Another is to train a 1-class support vector machine (SVM) and use it to classify samples. Here, we will compare against both of these standard techniques. To quantify performance, we will use the threshold-independent receiver operating characteristic (ROC) curves and area under curve (AUC) metric; this is a standard way of assessing binary classifier performance.  Both ROC and AUC are readily obtained using the MATLAB command {\tt perfcurve}. The kernel density estimates are obtained using the MATLAB command {\tt ksdensity} (and the KDE toolbox for higher dimensional examples \cite{kde_toolbox}). For the 1-class SVM, we use the MATLAB command {\tt fitcsvm} with the automatic tuning- and ``standardized data'' options enabled. 

We would like to emphasize to the reader, however, that our goal is not to estimate the probability density function. Our goal is simpler: to use the upper-bound on the probability of the observation to decide whether or not it is anomalous. To use $\mathbb{P}_S$ for anomaly detection we select $r>0$ to specify a neighborhood $S$ which contains an incoming measurement. The neighborhood can be a sphere or a box\footnote{$S$ in $\mathbb{P}_S$ can be any set of the form in (\ref{eq:semialgebraic_set}).} and is specified in the constraints  $f_j$ in $\mathbb{P}_S$. The $r$ can be selected as a function of the measurement noise, or otherwise small compared to the expected range of the data. In fact, $r$ can be set to zero if one wishes. In practice, we would like to use $r$ to account for measurement noise.

\subsection{Moment Estimation and Data Whitening}
A key advantage of the proposed method is that the information from a data set enters through the moment estimates which are computed using Equation (\ref{eq:raw_moment_estimate}). 

\begin{equation}
\gamma_\alpha \approx \frac{1}{N}\sum_{i=1}^N x_i^\alpha, \; \alpha\in\Gamma
\label{eq:raw_moment_estimate}
\end{equation}

It is known that it is difficult to estimate higher order moments since whatever noise is present in the data will be raised to higher powers as well. In view of Equation (\ref{eq:raw_moment_estimate}), one can see that if the data is poorly scaled or biased in some coordinate, $\gamma_\alpha$ will increase quickly in those directions and will therefore result in fewer available moments for those directions and may cause numerical problems for the SDP solver. Thus, it is sometimes useful to use the data whitening technique discussed in \cite{bishop_text} to transform the data so that it has zero mean and unit variance. This may be useful in lower-dimensional instances, where it is possible to use these higher order moments without the size of the moment matrices exceeding computational resources. 

Data whitening is a simple procedure that involves subtracting the mean of the data and multiplying by a transformation matrix, in the multivariate case, or dividing by the standard deviation in the univariate case. Incoming samples are then processed the same way, using the stored mean and transformation, before obtaining an upper-bound probability estimate.

\section{Examples}
In this section we present four examples\footnote{We used MATLAB R2013b/R2014a, CVX \cite{cvx}, and SeDuMi \cite{sedumi} to solve our examples.}. The first is a 1D example that shows that our approach compares favorably with traditional kernel probability density based anomaly detectors and with the 1-class SVM. We will demonstrate: that higher-order moments contain useful information and that this information helps best the other two approaches especially when there are fewer data points to learn from and when the dimensionality of the data increases.  These points will become apparent through the first three examples, which increase in dimension. The final example shows how our method can be used to recover contours from a noisy binary image.

\subsection{Example 1: A Bimodal Distribution} 
\begin{figure}[H]
\centering
\resizebox{\fig_size in}{!}{\includegraphics{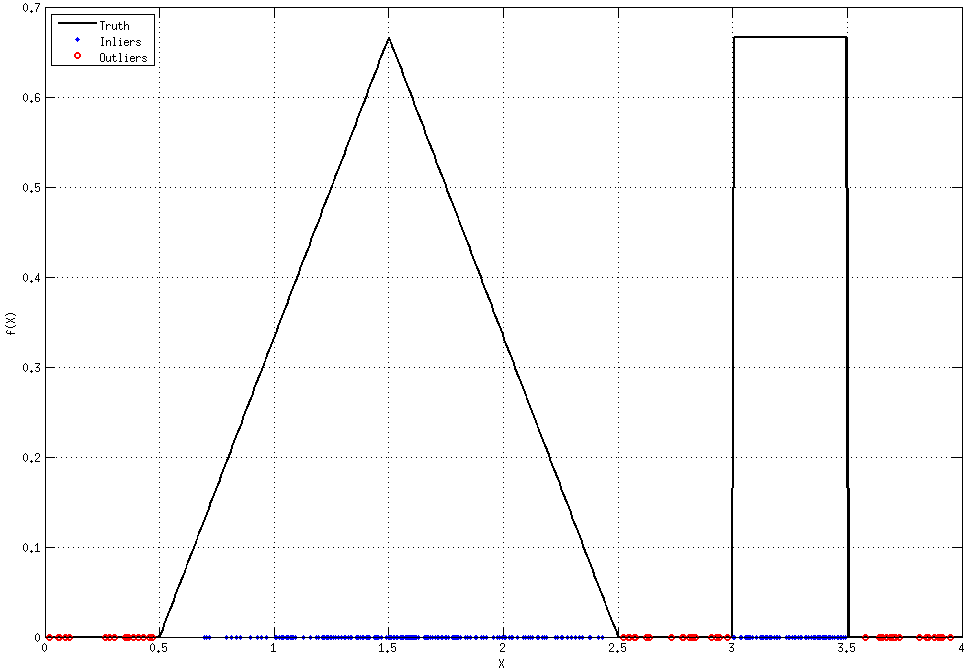}}
\caption{Bimodal distribution}
\label{fig:ex1_figure1}
\end{figure}

Consider the bimodal distribution \cite{duda_text} in Figure \ref{fig:ex1_figure1}. This distribution poses a challenge to density estimation tools like Parzen windows because the two modes are estimated best using different window sizes. 

Figure \ref{fig:ex1_figure1} shows the distribution and the test points, which consist of 300 inliers and 50 outliers. 

Table \ref{table:ex1_table} shows the AUC values for our experiments, which use $r=0.001$ for the moments classifier. The reader can see that  AUC generally increases as more moments are used and that the performance is better with fewer points. The difference, however, is small because this is a 1D example.
\begin{table}[H]
\centering
\caption{AUC Summary}
\begin{tabular}{ | c | c | c | c | c | c | c |}
\hline
\textbf{N} & \textbf{SVM} & \textbf{Parzen} & $M_2$ & $M_4$ & $M_6$ & $M_8$\\
\hline
50 & 0.9821 & 0.9686 & 0.8701 & 0.9715 & \textbf{0.9940} & 0.9929\\
\hline
100 & 0.9907 & 0.9807 & 0.8457 & 0.9882 & \textbf{0.9984} & 0.9957\\
\hline
300 & \textbf{0.9953} & 0.9781 & 0.8268 & 0.9783 & \textbf{0.9951} & 0.9945\\
\hline
\end{tabular}
\label{table:ex1_table}
\end{table}

\begin{figure}[H]
\centering
\resizebox{\fig_size in}{!}{\includegraphics{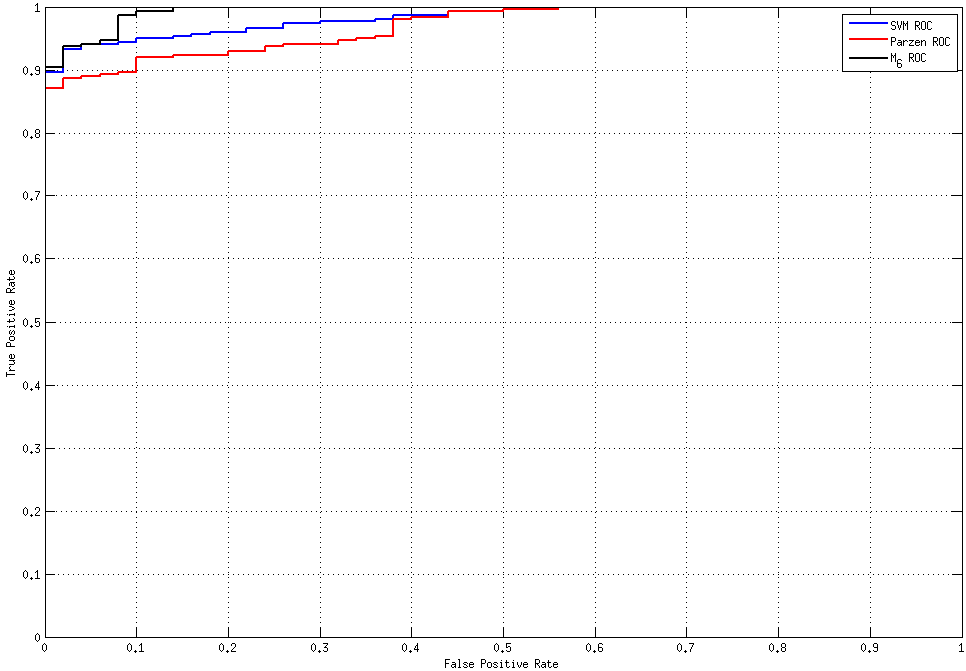}}
\caption{ROC curves when 50 points are used for training}
\label{fig:ex1_figure2}
\end{figure}

\begin{figure}[H]
\centering
\resizebox{\fig_size in}{!}{\includegraphics{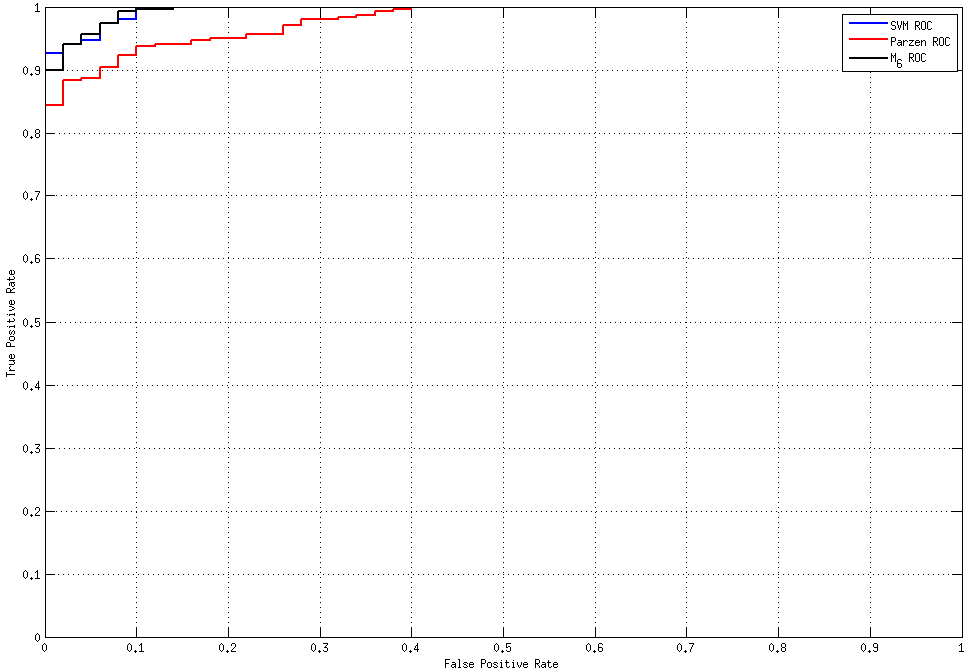}}
\caption{ROC curves when 300 points are used for training}
\label{fig:ex1_figure3}
\end{figure}

\subsection{Example 2: 2D Distribution} For higher dimensional datasets, the complexity of estimating PDFs increases greatly. Indeed, the arguably most popular engineering software tool, MATLAB, does not have a density estimation function for multi-dimensional datasets. 

Consider the distribution below 
\begin{align}
x_1 &= 2\cos(1.5 t) + \eta_1 \\
x_2 &= 4\sinc(0.9 t) + \eta_2 \notag
\end{align}
with $t{\in}[\pi/4,\pi]$, $\eta_1 {\sim} \mathcal{U}[-0.05,0.05]$, and $\eta_2 {\sim} \mathcal{U}[0,0.02]$. This distribution has a ``P'' shape and the closed portion presents a problem for the kernel and SVM approaches.

\begin{figure}[H]
\centering
\resizebox{\fig_size in}{!}{\includegraphics{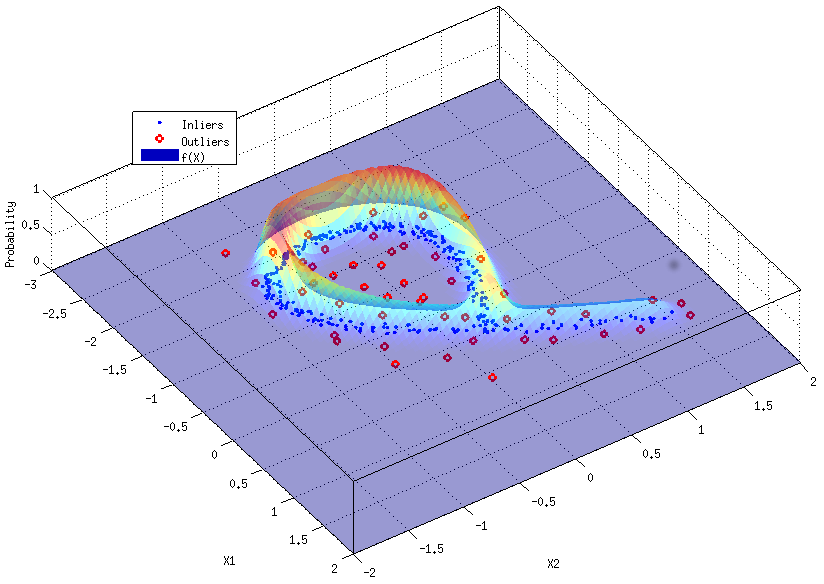}}
\caption{2D distribution and test points}
\label{fig:ex2_figure1}
\end{figure}

Figure \ref{fig:ex2_figure1} shows the probability surface and the 350 test points which include  50 outliers. 

Table \ref{table:ex2_table} shows the results of our experiments, using $r=0.001$ for the moments classifier. Again, the reader can see how the performance increases as more moment information is included except this time the difference in performance is greater because the data has another dimension. This trend will continue in the third example.

\begin{table}[H]
\centering
\caption{AUC Summary}
\begin{tabular}{ | c | c | c | c | c | c | c |}
\hline
\textbf{N} & \textbf{SVM} & \textbf{Parzen} & $M_2$ & $M_3$ & $M_4$ & $M_5$\\
\hline
100 & 0.9252 & 0.7843 & 0.6330 & 0.9722 & 0.9883 & \textbf{0.9916}\\
\hline
150 & 0.9437 & 0.8045 & 0.6575 & 0.9751 & 0.9903 & \textbf{0.9947}\\
\hline
300 & 0.9435 & 0.8335 & 0.6436 & 0.9705 & 0.9863 & \textbf{0.9933}\\
\hline
\end{tabular}
\label{table:ex2_table}
\end{table}

\begin{figure}[H]
\centering
\resizebox{\fig_size in}{!}{\includegraphics{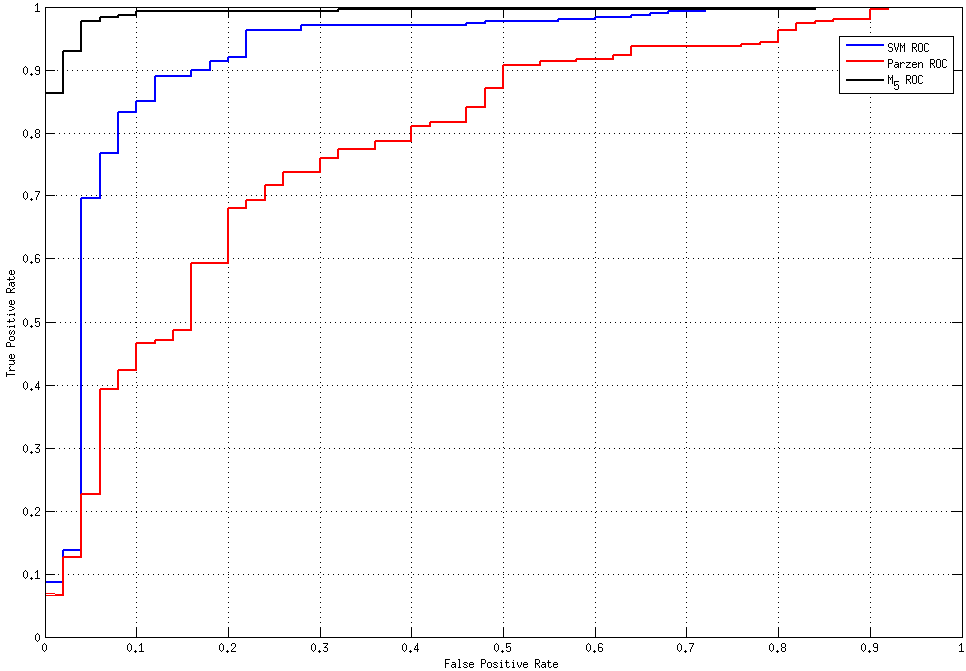}}
\caption{ROC curves when 100 points are used for training}
\label{fig:ex2_figure2}
\end{figure}

\begin{figure}[H]
\centering
\resizebox{\fig_size in}{!}{\includegraphics{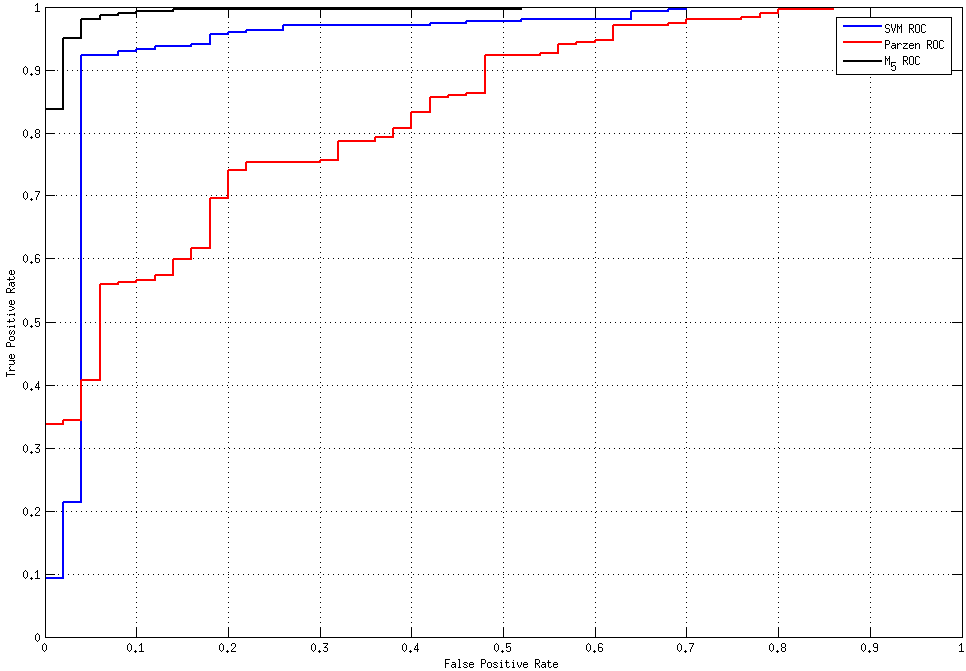}}
\caption{ROC curves when 300 points are used for training}
\label{fig:ex2_figure3}
\end{figure}

\subsection{Example 3: Swiss roll Distribution} 
\begin{figure}[H]
\centering
\resizebox{\fig_size in}{!}{\includegraphics{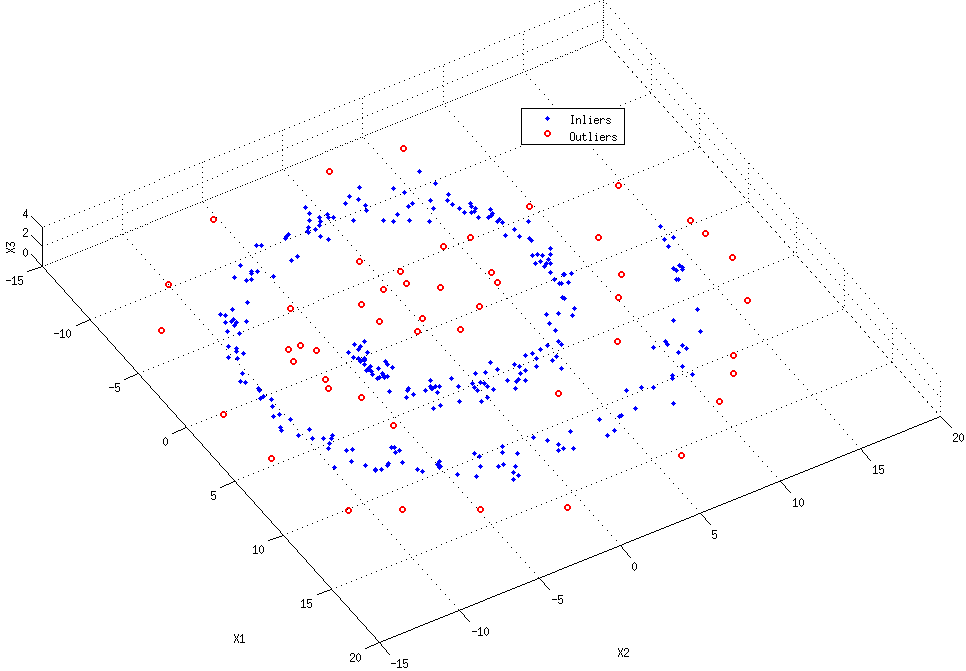}}
\caption{Swiss-roll distribution}
\label{fig:ex3_figure1}
\end{figure}

The following distribution is the well-known ``Swiss roll'' distribution. The following equations were used to create the test and training samples.

\begin{align}
x_1 &= y_1\cos (y_1) + \eta_1\\
x_2 &= y_1\sin (y_1) + \eta_2\\
x_3 &= y_2 + \eta_3
\end{align}
with $y_1\sim \mathcal{U}[\frac{3\pi}{2},\frac{9\pi}{2}]$, $y_2\sim \mathcal{U}[1,2]$, $\eta_1\sim\mathcal{N}(0,0.8^2)$, $\eta_2\sim\mathcal{N}(0,0.3^2)$, and $\eta_3\sim\mathcal{N}(0,0.25^2)$. Figure \ref{fig:ex3_figure1} shows the test points which consist of 300 inliers and 50 outliers.

Table \ref{table:ex3_table} shows the familiar trend of increasing AUC with additional moments and excellent performance even with just 100 training samples. This example illustrates how the other two standard techniques degrade as dimensionality increases. As with the previous examples, $r=0.001$ was used for the moments classifier and the recommended automatic tuning options were used for other techniques.

\begin{table}[H]
\centering
\caption{AUC Summary}
\begin{tabular}{ | c | c | c | c | c | c |}
\hline
\textbf{N} & \textbf{SVM} & \textbf{Parzen} & $M_2$ & $M_4$ & $M_6$\\
\hline
100 & 0.7887 & 0.8487 & 0.8099 & 0.9194 & \textbf{0.9228}\\
\hline
300 & 0.8001 & 0.9085 & 0.8201 & 0.9480 & \textbf{0.9903}\\
\hline
\end{tabular}
\label{table:ex3_table}
\end{table}

\begin{figure}[H]
\centering
\resizebox{\fig_size in}{!}{\includegraphics{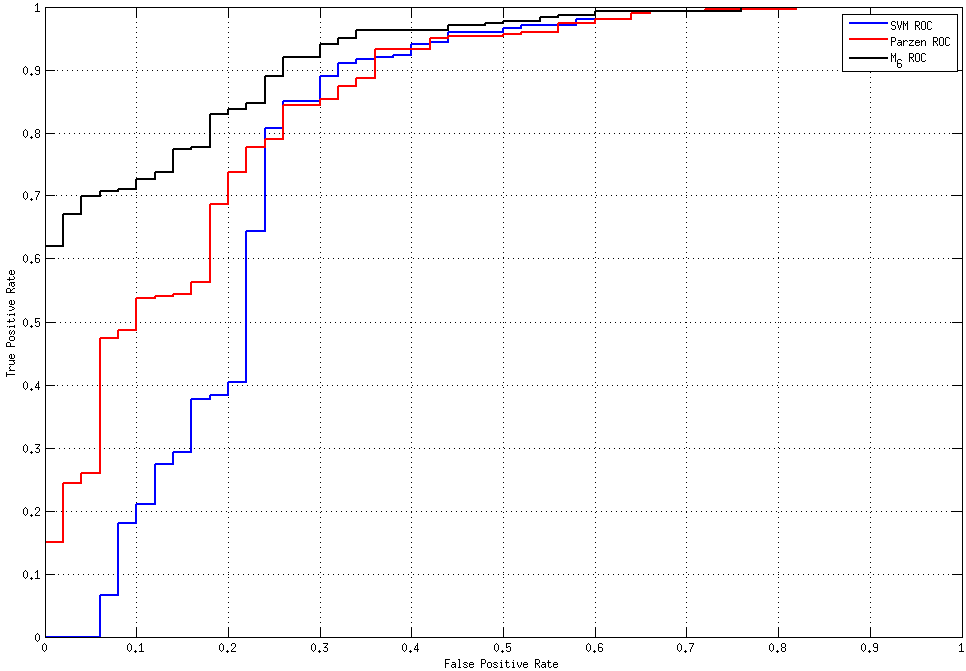}}
\caption{ROC curves when 100 points are used for training}
\label{fig:ex3_figure2}
\end{figure}

\begin{figure}[H]
\centering
\resizebox{\fig_size in}{!}{\includegraphics{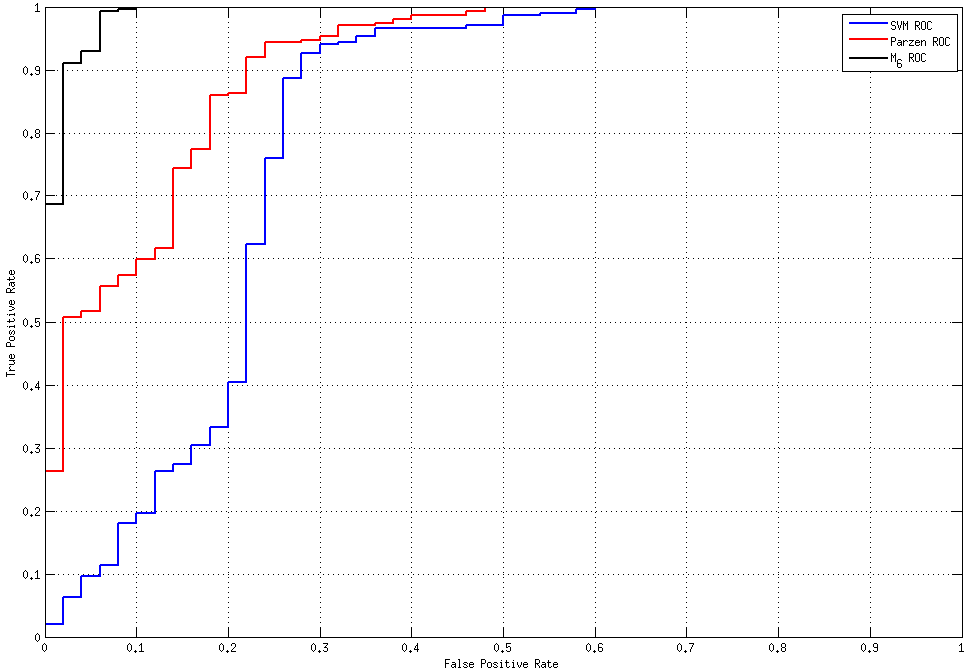}}
\caption{ROC curves when 300 points are used for training}
\label{fig:ex3_figure3}
\end{figure}

\subsection{Example 4: Contour Detection} 
\begin{figure}[H]
\centering
\resizebox{\fig_size in}{!}{\includegraphics{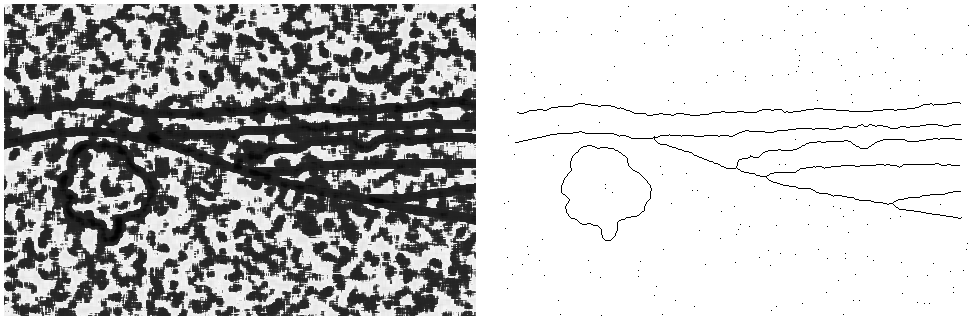}}
\caption{Upper-bound image and recovered contour}
\label{fig:prob_image_and_post_processed}
\end{figure}

The following example is included to show the potential of this method to be used in image processing applications. The contour image is from \cite{felzenszwalb_oberlin_arxiv,felzenszwalb_oberlin_nips}. In this example, we attempt to recover a contour corrupted by 5\% noise. We calculated the Haralick energy, in two directions and three scales, for (9x9) neighborhoods of each pixel of the clean image, computed $\gamma_\alpha$, and the upper-bound probability surface, $M_4$, using the noisy image with $r=0.05$. The three resulting probability images were then averaged. Since there are many more white pixels, the contour pixels should have a lower probability upper-bound. Further, as the Haralick energy is the sum of the squared elements of the co-occurrence matrix, computed for each 9x9 window, we also expect the noisy pixels to have higher probability than the contour pixels. The reader can see that this is the case, as shown in the left subfigure of Figure \ref{fig:prob_image_and_post_processed}. After inverting, thresholding, and post-processing with morphological operations, we obtained the contour shown on the right of Figure \ref{fig:prob_image_and_post_processed}, again, using just the information contained in $\gamma_\alpha$.
\begin{figure}[H]
\centering
\resizebox{\fig_size in}{!}{\includegraphics{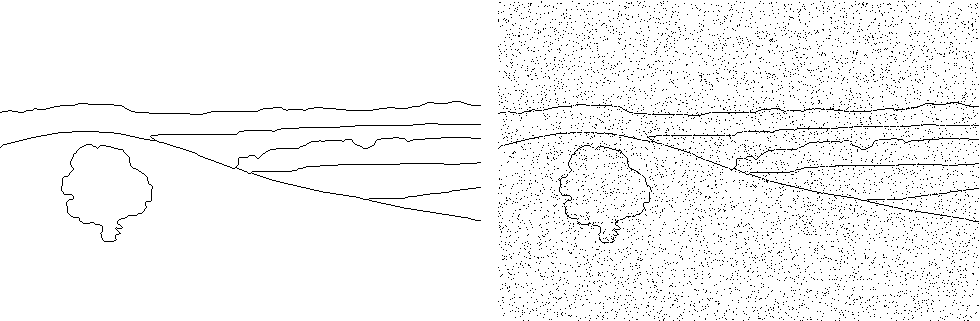}}
\caption{Clean and noisy images}
\label{fig:clean_and_noisy}
\end{figure}

\section{Conclusions}
This paper illustrates the ability  of very recent developments in polynomial optimization and generalized moment problems to provide 
an elegant, computationally tractable, solution to anomaly detection problems. The proposed approach uses information from a set of moments of a probability distribution to directly provide an upper-bound $\rho$ on the probability of observing a particular sample (modulo measurement noise), without explicitly modeling the distribution. A salient feature of this technique is that all the data, whether obtained using 1 or 1000 observations, enter the problem through a finite sequence of moments, that can be updated as more data is collected, without affecting the computational complexity (since new observations affect the value, but not the size of the moment matrices). Furthermore, the detectors obtained  by solving $\mathbb{P}_S$  reached peak performance with fewer samples than the kernel and SVM methods. 

Our examples implicitly assume that the observations contain enough information to make good decisions, that is, we have assumed that the samples are \emph{good} features. When anomalous observations are available, it may be possible to construct optimal features to separate classes -- in the same way support (or relevance) vector machine methods find separating curves -- by maximizing the distance between the moment vectors of the two classes. We think this is an interesting direction for future research.

\bibliographystyle{plain}
\bibliography{references}

\end{document}